\documentclass[12 pt,twoside]{article}
\usepackage{graphicx}

\newtheorem{theorem}{\bf Theorem}[section]

\newtheorem{lemma}[theorem]{\bf Lemma}

\newtheorem{emp}[theorem]{\bf Claim}

{{\footnotesize\sc }}
\input{amssym}

\begin{document}
\vspace{2 cm}
\title{{\large Tur\'{a}n numbers of complete $3$-uniform Berge-hypergraphs}}
\vspace{4cm}
\bigskip\author{\small L. Maherani$^{\textrm{a}}$,  M. Shahsiah$^{\textrm{b,c}}$ \\
\footnotesize  $^{\textrm{a}}$ Department of Mathematical Sciences,
Isfahan
University of Technology,\\ \footnotesize Isfahan, 84156-83111, Iran\\
{\small $^{\textrm{b}}$Department of Mathematics, Alzahra University,}\\
{\small P.O. Box 1993891176, Tehran, Iran}\\
{\small $^{\textrm{c}}$School of Mathematics, Institute for Research in Fundamental Sciences (IPM),}\\
{\small P.O. Box 19395-5746, Tehran, Iran }\\
\footnotesize {l.maherani@math.iut.ac.ir,
 shahsiah@ipm.ir}}
\date {}

\footnotesize\maketitle

\begin{abstract}\rm{}
\medskip
\footnotesize
Given a family $\mathcal{F}$ of $r$-graphs, the Tur\'{a}n number of $\mathcal{F}$ for a given
positive integer $N$, denoted by $ex(N,\mathcal{F})$, is the maximum number of edges of an $r$-graph on $N$ vertices
that does not contain any member of $\mathcal{F}$ as a subgraph.  For given $r\geq 3$, a  complete $r$-uniform Berge-hypergraph, denoted by { ${K}_n^{(r)}$}, is an $r$-uniform hypergraph of order $n$  with the core sequence $v_{1}, v_{2}, \ldots ,v_{n}$ as the vertices and  distinct edges $e_{ij},$ $1\leq i<j\leq n,$ where every $e_{ij}$ contains both   $v_{i}$ and $v_{j}$. Let   $\mathcal{F}^{(r)}_n$ be the family of complete $r$-uniform Berge-hypergraphs of order $n.$  We determine precisely $ex(N,\mathcal{F}^{(3)}_{n})$ for $n \geq 13$. We also find the extremal hypergraphs avoiding $\mathcal{F}^{(3)}_{n}$.
\vspace{.5cm}
 \\{ {Keywords}:{ \footnotesize Tur\'{a}n number, Extremal hypergraph, Berge-hypergraph.\medskip}}
\noindent
\\{\footnotesize {AMS Subject Classification}:  05C65, 05C35, 05D05.}
\end{abstract}
\small
\medskip
\section{\normalsize{Introduction}}

A {\it hypergraph} $\mathcal{H}$ is a pair $\mathcal{H}=(V,E)$, where $V$ is a finite non-empty set (the set of
vertices) and $E$ is a collection of distinct non-empty subsets of $V$ (the set of edges). We denote by $e(\mathcal{H})$ the number of edges of $\mathcal{H}.$
 An {\it $r$-uniform
hypergraph} or {\it $r$-graph} is a hypergraph such that all its edges have size $r$.   A {\it complete $r$-uniform  hypergraph} of
order $N$, denoted by { $\mathcal{K}_N^r$}, is a hypergraph   consisting of
all the $r$-subsets of a set $V$ of cardinality $N$. For a family $\mathcal{F}$ of $r$-graphs, we say that the hypergraph $\mathcal{H}$ is $\mathcal{F}$-free if $\mathcal{H}$  does not contain any member of $\mathcal{F}$ as a subgraph.
 Given a family $\mathcal{F}$ of $r$-graphs, the {\it Tur\'{a}n number} of $\mathcal{F}$ for a given
positive integer $N$, denoted by $ex(N,\mathcal{F})$, is the maximum number of edges of an  $\mathcal{F}$-free  $r$-graph on $N$ vertices.
 An $\mathcal{F}$-free $r$-graph $\mathcal{H}$ on  $N$ vertices is {\it extremal hypergraph} for $\mathcal{F}$ if $e(\mathcal{H})= ex(N,\mathcal{F})$. These are  natural generalizations of the classical Tur\'{a}n number for
$2$-graphs \cite{turan}.
 For given $n, r \geq 2$,  let $\mathcal{H}^{(r)}_n$ be the family of $r$-graphs $F$ that have at most ${n \choose 2}$ edges, and have some set $T$ of size $n$ such that every pair of vertices in $T$ is contained in some edge of $F$.
Let  the $r$-graph  $H^{(r)}_n \in \mathcal{H}^{(r)}_n$ be obtained from the complete $2$-graph $\mathcal{K}_{n}^{2}$ by enlarging each edge with a new set of $r-2$ vertices. Thus $H^{(r)}_n$ has $(r-2){n \choose 2}+n$ vertices and ${n \choose 2}$ edges. For given $n \geq 5$ and  $r\geq 3$, a {\it complete $r$-uniform Berge-hypergraph} of order $n$, denoted by { ${K}_n^{(r)}$}, is an $r$-uniform hypergraph with the core sequence $v_{1}, v_{2}, \ldots ,v_{n}$ as the vertices and ${n \choose 2}$ distinct edges $e_{ij},$ $1\leq i<j\leq n,$ where every $e_{ij}$ contains both   $v_{i}$ and $v_{j}$. Note that a complete $r$-uniform Berge-hypergraph is not determined uniquely as there are no constraints on how the $e_{ij}$'s intersect outside $\{v_{1}, v_{2}, \ldots ,v_{n}\}$.
\\

Extremal graph theory is that area of combinatorics
which is concerned with finding the largest, smallest, or otherwise optimal
structures with a given property. There is a long history in
the study of extremal problems concerning hypergraphs. The first such result is due to Erd\H{o}s, Ko and Rado \cite{Erdos-ko-rado}.\\

In contrast to the graph case, there are comparatively few known results  on
the hypergraph Tur\'{a}n problems. In the
 paper in which Tur\'{a}n proved his classical theorem on the extremal numbers
for complete graphs \cite{turan}, he posed the natural question of determining the Tur\'{a}n number of the complete $r$-uniform  hypergraphs.  Surprisingly, this
problem remains open in all cases for $r > 2$, even up to asymptotics. Despite the lack of progress on the Tur\'{a}n problem for dense hypergraphs, there are considerable results on certain sparse hypergraphs.
Recently, some interesting results were
obtained on the exact value of extremal number of paths and cycles in hypergraphs.
F{\"u}redi et al.  \cite{loose pathI} determined the extremal number of $r$-uniform loose paths of length $n$ for $r\geq 4$ and large $N$.
 They also  conjectured a similar result for $r = 3$.  F{\"u}redi and Jiang \cite{loose cycleI} determined the extremal function of loose cycles  of length $n$ for  $r\geq 5$ and large $N.$ Recently, Kostochka et al. \cite{loose paths and cycles} extended these results to $r=3$ for loose paths and $r=3,4$ for loose cycles. Gy{\H{o}}ri et al. \cite{Berge paths} found the extremal numbers of $r$-uniform hypergraphs avoiding Berge paths of length $n$. Their results substantially extend earlier results of Erd\H{o}s  and Gallai \cite{Erdos-Gallai} on extremal number of paths in graphs. Let $\mathcal{C}_{n}^{(r)}$
denote the family of $r$-graphs that are Berge cycles of length $n$.  Gy{\H{o}}ri and Lemons \cite{Berge cyclesI, Berge cycles2}
showed that for all $r \geq 3$ and $ n \geq 3$, there exists a positive constant $c_{r,n},$ depending on $r$ and $n,$ such that
$$ex(N,\mathcal{C}_{n}^{(r)}) \leq c_{r,n} N^{1+\frac{1}{\lfloor \frac{n}{2} \rfloor}}.$$
\vspace{0.2 cm}

Let $N$, $n$, $r$ be integers, where $N\geq n >r$ and $r \geq 2$. Also let $ T_r(N,n-1)$ be the complete $r$-uniform $(n-1)$-partite hypergraph with $N$ vertices and $n-1$  parts $V_1,V_2,...,V_{n-1}$ whose partition sets differ in size by at most 1. Suppose that  $t_r(N,n-1)$ denotes the number of edges of $T_r(N,n-1)$. If $N=\ell(n-1)+j$, where $\ell \geq 1$ and $1 \leq j\leq n-1$, then  it is straightforward to see that
$$t_r(N,n-1)= \sum _{i=0}^{r} \ell ^{r-i} {j \choose i}{n-1-i \choose r-i}.$$

In 2006, Mubayi \cite{mubay} showed that the unique largest $\mathcal{H}_{n}^{(r)}$-free $r$-graph on $N$ vertices is $T_r(N,n-1)$.
 Settling a conjecture  of Mubayi  in \cite{mubay}, Pikhurko \cite{pikh} proved that  there exists $N_0$ so that  the
Tur\'{a}n numbers of ${H}_{n}^{(r)}$ and $\mathcal{H}_{n}^{(r)}$ coincide for all  $N>N_0.$
Let   $\mathcal{F}^{(r)}_n$ be the family of complete $r$-uniform Berge-hypergraphs of order $n.$ Because ${H}_{n}^{(3)} \in \mathcal{F}_{n}^{(3)}$, the  Pikhurko's result  \cite{pikh}  implies that $ex (N,\mathcal{F}_{n}^{(3)}) \leq t_3(N,n-1)$ for sufficiently large $N$.
In this paper, for $N \geq 13$,  we show  that $ex (N,\mathcal{F}_{n}^{(3)})=t_3(N,n-1)$ and $ T_3(N,n-1)$  is the unique extremal hypergraph for $\mathcal{F}_{n}^{(3)}$. More precisely, we prove the following theorem.\\

\begin{theorem}\label{main}
Let $N,n$ be integers so that $N\geq n\geq 13$. Then
$$ex(N,\mathcal{F}_n^{(3)})=t_3(N,n-1).$$
Furthermore, the unique extremal hypergraph for $\mathcal{F}_n^{(3)}$ is $T_3(N,n-1)$.
\end{theorem}
\vspace{0.5 cm}

%

First we show that $ex (N,\mathcal{F}_{n}^{(r)}) \geq t_r(N,n-1)$. To see that, consider an arbitrary sequence $v_1,v_2,...,v_n$ of the vertices of $ T_r(N,n-1)$. By the pigeonhole principle, there exists some part $V_h$, $1 \leq h \leq n-1$, in $T_r(N,n-1)$ containing at least two vertices of this sequence. Since every edge of $ T_r(N,n-1)$ includes at most one vertex of each part $V_i$, $1 \leq i \leq n-1$,  This sequence can not be the core sequence of a $K_{n}^{(r)}$. Hence $T_r(N,n-1)$ is $\mathcal{F}_{n}^{(r)}$-free and

\begin{equation}\label{lbound}
ex (N,\mathcal{F}_{n}^{(r)}) \geq t_r(N,n-1), \ \ \ \ \ \ r\geq 3.
\end{equation}

 Therefore, in order to clarify Theorem \ref{main}, it  suffices to show that $ex (N,\mathcal{F}_{n}^{(3)}) \leq t_3(N,n-1)$ and $T_3(N,n-1)$ is the only $\mathcal{F}_{n}^{(3)}$-free hypergraph with $N$ vertices and $t_3(N,n-1)$ edges. Here, we give a proof by induction on the number of vertices. More precisely, we prove Theorem \ref{main} in three steps. First, we show that Theorem \ref{main} holds for $N=n$ (see Theorem \ref{n}). Then, in Theorem \ref{l=1}, we demonstrate that it is true for $n\leq N \leq 2n-2.$ Finally, using Theorem \ref{n} and Theorem \ref{l=1}, we show that the desired holds for all $N\geq n$ (Section  3).\\

 \noindent
\textbf{Conventions and Notations:}
For an $r$-uniform hypergraph $\mathcal{H}=(V,E)$, the complement hypergraph of $\mathcal{H}$, denoted by $\mathcal{H}^c$, is the hypergraph on $V$ so that $E(\mathcal{H}^c)= {V \choose r}\setminus E$. Also we say that $X\subseteq V$ is an independent set of $\mathcal{H}$ if for any pair  $ v,v' \in X$, there is no edges in $E$ containing both of $v$ and $v'$. For $U\subseteq V$ we denote by $\mathcal{H}[U]$ the subgraph of $\mathcal{H}$ induced by the edges of $U$. For $U, W \subseteq V$, The hypergraph $\mathcal{H}[U,W]$ is the subgraph of $\mathcal{H}$ induced by the edges of $\mathcal{H}$ intersecting both $U$  and $W.$ For  a vertex $v \in V$, the degree of $v$ in $\mathcal{H}$, denoted by $d_\mathcal{H}(v)$, is the number of edges in $\mathcal{H}$ containing $v$.  Also $\mathcal{H}-v$ is the subhypergraph of $\mathcal{H}$ obtained by deleting of $v$ and all the edges containing it.\\

\section{Preliminaries }
 In this section, we present some
results that will be used in the follow up section. Let $\mathcal{A}=\{A_{1},A_{2},...,A_{n}\}$ be a family of subsets of a set $X$. A system of distinct representatives,
or SDR, for the family $\mathcal{A}$, is a set $\{a_1,a_2,...,a_n\}$ of elements of $X$ satisfying
two following conditions:
\begin{itemize}
\item[$\bullet$] $a_i \in A_i$ \ \ \ \ \ \ $i=1,...,n$,
\item[$\bullet$] $a_i \neq a_j$  \ \ \ \ \ \ $i \neq j$.
\end{itemize}


\medskip
\begin{lemma}\label{sdr}
Let $U=\{u_1,u_2,...,u_m\}$, $m\geq 5$ and $x \notin U$. Also, let $\mathcal{A}=\{A_{1},A_{2},...,A_{m}\}$ be a family of sets so that $\vert A_{1}\vert \leq \vert A_{2}\vert \leq ... \vert A_{m}\vert $ and
$A_{i} \subseteq \{ B :\  B=\{x, u_i ,u_k \},\ k\neq i \}$  for $ 1\leq i \leq m.$
If $\mathcal{A}$ has no SDR, then $$\vert \bigcup _{i=1}^{m} A_{i} \vert \leq {m-1 \choose 2}$$ and equality holds if and only if $A_{1}= \emptyset$ and
$$A_{i} = \{ B :\  B=\{x, u_i ,u_k \},\ k\neq 1,i \}  \ \ \ \ \ \ 2\leq i \leq m.$$
\end{lemma}
\noindent\textbf{Proof. }  Since  $\mathcal{A}$ contains no SDR, using the Hall's theorem \cite{hall}, for some  $q$, $1\leq q\leq m$, we have $\vert \bigcup _{i=1}^{q} A_{i} \vert \leq q-1$. So $\vert \bigcup _{i=1}^{m} A_{i} \vert \leq f(q),$ where $f(k)= k-1 + {m-k \choose 2}$, for $1 \leq k \leq m$. On the other hand, one can easily see that $f(1)>f(k)$, for $2\leq k \leq m$. Therefore
$$\vert \bigcup _{i=1}^{m} A_{i} \vert \leq f(1)={m-1 \choose 2}$$ and the equality holds if and only if $A_{1}= \emptyset$ and
$$A_{i} = \{ B :\  B=\{x, u_i ,u_k \},\ k\neq 1,i \} \ \ \ \ \ \ 2\leq i \leq m.$$
$\hfill\blacksquare$\\

In order to state our main results we need some definitions.  Let  $\mathcal{H}=(V,E)$ be an $r$-uniform hypergraph, where $V=\{v_1,v_2,...,v_n\}$ and $E=\{e_1,e_2,...,e_m\}$.
We denote by $B(\mathcal{H}),$ the bipartite graph with parts $X$ and $Y$
  so that  $X=\{v_iv_k :\ i< k\ \  {\rm and} \  \  v_i,v_k \in V(\mathcal{H})\}$, $Y=E(\mathcal{H})$ and $v_iv_k$ is adjacent to $e_h$ if and only if $\{v_i,v_k \}\subseteq e_h$, for every $v_iv_k\in X$ and $e_h \in Y$. For every $v_iv_k\in X$, $d_{B(\mathcal{H})}(v_iv_k)$ is the number of edges in $B(\mathcal{H})$ containing  $v_i v_k$. 
A matching of $X$ in $B(\mathcal{H})$ is matching that saturates all vertices of $X$.  Note that, every matching of $X$  in $B(\mathcal{H})$ is equivalent to a complete $r$-uniform Berge-hypergraph with core sequence $v_1,v_2,...,v_n$. \\

\noindent Now, we demonstrate that Theorem \ref{main} holds for $N=n$.
\begin{theorem}\label{n}
Let $n\geq 13$ be an integer. The hypergraph $T_3(n,n-1)$ is the only  $\mathcal{F}_{n}^{(3)}$-free hypergraph with $n$ vertices and $ex(n,\mathcal{F}_{n}^{(3)})$ edges.
\end{theorem}

\noindent\textbf{Proof. }Assume that $\mathcal{H}$ is an $\mathcal{F}_{n}^{(3)}$-free hypergraph with $n$ vertices and $ex(n,\mathcal{F}_{n}^{(3)})$ edges. Let $V(\mathcal{H})=\{v_1,v_2,...,v_n\}$. First, suppose that there is a vertex $v\in V(\mathcal{H})$, say $v_n$, so that $d_{\mathcal{H}}(v_n) \leq {n-2 \choose 2}$. Therefore
\begin{equation}\label{up1}
  e(\mathcal{H})  =d_{\mathcal{H}}(v_n) +e(\mathcal{H}-v_n) \leq  {n-2 \choose 2} +{n-1 \choose 3}=t_3(n,n-1).
\end{equation}

\noindent So by (\ref {lbound}) and (\ref{up1}), we have
$$ex(n,\mathcal{F}_{n}^{(3)})  = t_3(n,n-1).$$
Therefore  $d_{\mathcal{H}}(v_n) ={n-2 \choose 2}$ and $e( \mathcal{H}-v_n )={n-1 \choose 3}$. So $\mathcal{H}-v_n \cong \mathcal{K}_{n-1}^3$ and clearly there is a copy of   $K_{n-1}^{(3)}$ with the core sequence $v_1,v_2,...,v_{n-1}$ in  $\mathcal{H}-v_n$.
Set $x=v_n$, $U=\{v_1,v_2,...,v_{n-1}\}$ and  $\mathcal{A}=\{A_{1},A_{2},...,A_{n-1}\}$, where $$A_{i} = \{ B :\  B\in E(\mathcal{H}),\  \{x,v_{i}\} \subseteq B\} \ \ \ \ \ \  1\leq i\leq n-1.$$
Note that $d_{\mathcal{H}}(v_n)=\vert \bigcup _{i=1}^{n-1} A_{i} \vert= {n-2 \choose 2}.$
Since $ \mathcal{H}$ is $\mathcal{F}_{n}^{(3)}$-free and there is a copy of   $K_{n-1}^{(3)}$ in $ \mathcal{H}-v_n$, $\mathcal{A}$ has no SDR.
Now,  using Lemma \ref{sdr}, we have $\mathcal{H}\cong T_3(n,n-1)$.\\*

 Now suppose that  for every vertex $v\in V(\mathcal{H})$, $d_{\mathcal{H}}(v) \geq {n-2 \choose 2}+1$. Set $G=B(\mathcal{H})$. So we may assume that $G=[X,Y]$, where
 $$X=\{u_{ik}=v_iv_k :\ i< k\ \  {\rm and} \  \  v_i,v_k \in V(\mathcal{H})\}$$ and $Y=E(\mathcal{H})$. Since, by (\ref{lbound}), $|Y|\geq  {n-1 \choose 3} +{n-2 \choose 2},$ we have $\vert X\vert \leq \vert Y\vert$. Let $X=X_1 \cup X_2$, where $X_{1}=\{u \in X :\  d_{G}(u)\leq 4\}$ and  $X_2 =X\setminus X_1$. 
Recall that every matching of $X$ in $G$ is equivalent to a $K_{n}^{(3)}$ in $\mathcal{H}$.
We have two following cases.\\

\noindent{\bf Case 1.} $X_1 =\emptyset$.\\
Since for every $y\in Y$ and $u\in X$, we have $d_{G}(y)=3$ and $d_{G}(u)\geq 5$, the Hall's theorem \cite{hall}  guarantees the existence of a matching of $X$, a contradiction.\\\\

\noindent{\bf Case 2.} $X_1 \neq \emptyset$.\\
Let $X_1 =  \{v_{i_1}v_{i'_1},v_{i_2}v_{i'_2},...,v_{i_t}v_{i'_t}\}$. We show that the following claim holds.

\begin{emp}\label{pairdisj1}
The elements of $X_1$ are pairwise disjoint.
\end{emp}
\noindent\textbf{Proof of Claim \ref{pairdisj1}}. Suppose to contrary that for $2 \leq s \leq t$, $ \{wv_{i'_1},wv_{i'_2},...,wv_{i'_s}\} \subseteq X_1$. So $ d_{\mathcal{H}}(w) \leq f(s),$ where $f(k)= 4k + {n-k-1 \choose 2}$ is a function on $k$,  $2 \leq k \leq t \leq n-1$.  Using $n\geq 13$, it is straightforward to see that the absolute maximum of $f(k)$ occurs in point $k=2$. Hence $$ d_{\mathcal{H}}(w) \leq f(2) = 8 + {n-3 \choose 2}.$$
Since $ 8 + {n-3 \choose 2} < {n-2 \choose 2}+1$ for $n\geq 13$, we have $ d_{\mathcal{H}}(w)  <  {n-2 \choose 2}+1$. That is a contradiction to our assumption.
$\hfill\square$\\

\noindent Since for every vertex  $v \in V(\mathcal{H})$, we have $ d_{\mathcal{H}}(v)  \geq  {n-2 \choose 2}+1$, so for any two vertices $x,y \in V(\mathcal{H})$, there is at least one edge in $E(\mathcal{H})$ containing both of $x$ and $y$.
 So $d_{G}(v_{i_l}v_{i'_l}) \geq 1$ for every $1 \leq l \leq t$. On the other hand, by  Claim \ref{pairdisj1}, the elements of $X_1$ are pairwise disjoint. Therefore  $G$ contains a matching $M_1$ of $X_1$. Suppose that $G'=[X_2,Y']$ is the subgraph of $G$ so that $Y' \subset Y$ is obtained by deleting the vertices of $M_1$. Note that for every $u \in X_2$ and $y \in Y'$, we have  $d_{G'}(u) \geq 3$ and  $d_{G'}(y) \leq 3$. Therefore the Hall's theorem  \cite{hall} implies the existence of a matching $M_2$ of $X_2$ in $G'$. This is a contradiction, since $M_1 \cup M_2$ is a matching of $X$ in $G$. This contradiction completes the proof.
$\hfill\blacksquare$

\begin{theorem}\label{l=1}
Let $n\geq 13$ and $N,n$ be integers so that $n\leq N \leq 2n-2$.
 Also, let $\mathcal{H}$ be an  $\mathcal{F}_{n}^{(3)}$-free hypergraph with $N$ vertices and $ex(N,\mathcal{F}_{n}^{(3)})$ edges. Then $ e(\mathcal{H}) = t_3(N,n-1)$ and $\mathcal{H}\cong T_3(N,n-1)$.
\end{theorem}
\noindent\textbf{Proof. } Let $N=n-1+j$, where $1\leq j \leq n-1$. We apply induction on $j$. Using Theorem \ref{n}, the basic step $j=1$ is true. For the induction step, let $j >1$. Set
$$d={n-2 \choose 2}+(j-1)(n-3)+{j-1 \choose 2}.$$
First suppose that there is a vertex $x \in V(\mathcal{H})$ so that $ d_{\mathcal{H}}(x)  \leq d$. So using the induction hypothesis, we have
$$e( \mathcal{H}) = d_{\mathcal{H}}(x)  + e( \mathcal{H}-x) \leq d+t_3(N-1,n-1) =t_3(N,n-1).$$
Therefore  by (\ref{lbound}), we conclude that  $e(N,\mathcal{F}_{n}^{(3)}) =t_3(N,n-1)$. Hence $d_{\mathcal{H}}(x) =d$ and  $e( \mathcal{H}-x) = t_3(N-1,n-1)$. So, using the induction hypothesis,  $\mathcal{H}-x \cong T_3(N-1,n-1)$. Hence we may assume that  $\mathcal{H}-x$ is a complete  $3$-uniform $(n-1)$-partite hypergraph  with parts $V_1,V_2,...,V_{n-1}$, where
$$
 V_i = \left\lbrace
\begin{array}{ll}
\{v_i,x_i\}  & \ \ \ \  1\leq i \leq j-1,\vspace{.5 cm}\\
\{v_i\}  & \ \ \ \   j\leq i \leq n-1.
\end{array}
\right.\vspace{.2 cm}
$$
Let $\mathcal{H}'$ be the induced subgraph of $\mathcal{H}-x$ on $\{v_1,v_2,...,v_{n-1}\}.$ According to the construction of $\mathcal{H}-x$, we have  $\mathcal{H}'\cong \mathcal{K}_{n-1}^3$ and so there is a copy of 
 $K_{n-1}^{(3)}$ with core sequence $v_1,v_2,...,v_{n-1}$  in $ \mathcal{H}'$. Set $U=\{v_1,v_2,...,v_{n-1}\}$ and $\mathcal{A}=\{A_{1},A_{2},...,A_{n-1}\}$, where for $1\leq i\leq n-1$,
$$A_{i} = \{ e \in E(\mathcal{H}) : \  e=\{x,v_{i},v_k\}, \ \ k\neq i\}.$$

\noindent For a vertex $v \in V(\mathcal{H})$, we denote by $E_v$ the set of  edges of $\mathcal{H}$ containing $v$. Clearly we have
\begin{equation}\label{dx}
d_{\mathcal{H}}(x)= |E_x| = |E_1|+ |E_2|+ |E_3|,
\end{equation}
where
$$E_i= \{e \in E_x:\ \  \vert e \cap \{x_1,x_2,...,x_{j-1}\}\vert =i-1\}, \ \ \ \ \ \ 1\leq i \leq 3.$$

\noindent We have the following claim.
\begin{emp}\label{Ex}
\noindent
\begin{itemize}
\item[{\rm (i)}] $\vert E_1\vert \leq {n-2 \choose 2}.$
\item[{\rm (ii)}] $ \vert E_2\vert \leq (j-1)(n-3).$
\item[{\rm(iii)}] $\vert E_3\vert \leq {j-1 \choose 2}.$
\end{itemize}
\end{emp}
\noindent\textbf{Proof of Claim \ref{Ex}}. (i)
Clearly $\vert E_1\vert =\vert \bigcup _{i=1}^{n-1} A_{i} \vert$. If $\mathcal{A}$ contains an SDR, then $x,v_1,v_2,...,v_{n-1}$ is the core sequence of a copy of  $K_n^{(3)}$ in $\mathcal{H}$, a contradiction. So, using Lemma \ref{sdr},
$$\vert E_1\vert =\vert \bigcup _{i=1}^{n-1} A_{i} \vert \leq {n-2 \choose 2}.$$

\noindent (ii) For $1 \leq k \leq j-1$, set $$B_k=\{e \in E_2:\ \{x,x_k\}\subseteq  e\}.$$
We demonstrate that  for $1 \leq k \leq j-1$, $\vert B_k \vert \leq n-3$ and so
$$\vert E_2\vert =\vert \bigcup _{k=1}^{j-1} B_k \vert \leq (j-1)(n-3).$$
Because of the similarity, it suffices to show that $\vert B_1 \vert \leq n-3$. Suppose not. So $\vert B_1\vert \geq n-2.$ On the other hand, the construction of $\mathcal{H}-x$ and the fact that $\mathcal{F}_{n}^{(3)}\nsubseteq \mathcal{H}$ imply that every edge in $E_x$ contains at most one vertex of each $V_i,$ for $1\leq i \leq n-1.$ Hence $|B_1|=n-2$ and 
$$B_1=\{\{x,x_1,v_2\}, \{x,x_1,v_3\},...,\{x,x_1,v_{n-1}\}\}.$$
In this case, there is no edge in $E(\mathcal{H})\setminus B_1$ containing  both of $x$ and $v_i,$ for $2 \leq i \leq n-1$. To see it, suppose that  $f=\{x,v_2,u\}\in E(\mathcal{H})\setminus B_1$. Let $\mathcal{H}''$ be the induced subgraph of $\mathcal{H}-x$ on $\{x_1,v_2,...,v_{n-1}\}.$ By the construction of $\mathcal{H}-x$, we have  $\mathcal{H}''\cong \mathcal{K}_{n-1}^3$ and so
 $\mathcal{H}''$ contains a $K_{n-1}^{(3)}$, say $\mathcal{K}'$.  Hence $x, x_1,v_2,...,v_{n-1}$ represents the core sequence of a $K_{n}^{(3)}$ in $\mathcal{H}$  with the following edge assignments. Set $e_{xx_1}=\{x,x_1,v_2\}$, $e_{xv_2}=f$, $e_{xv_i}=\{x,x_1,v_i\}$ for $3 \leq i \leq n-1$ and  other edges are selected from $E(\mathcal{K}')$. That is a contradiction to our assumption. Therefore the set of edges in $\mathcal{H}$ containing $x$ and $v_1$ is a subset of the following set:
$$S=\{\{x,v_1,x_2\}, \{x,v_1,x_3\},...,\{x,v_1,x_{j-1}\}\}.$$ Hence
$$d_{\mathcal{H}}(x) \leq |B_1|+ |S|+{j-1 \choose 2}=
(n-2)+(j-2)+{j-1 \choose 2} < d.$$ This contradiction demonstrates that $\vert B_1 \vert \leq n-3$ and so $\vert E_2\vert \leq (j-1)(n-3)$.\\

\noindent (iii) This case is trivial.
$\hfill\square$\\

\noindent Since $d_{\mathcal{H}}(x)=d$, using (\ref{dx}) and Claim \ref{Ex}, we have
\begin{equation}\label{ddd}
 \vert E_1\vert = {n-2 \choose 2},\ \ \ \vert E_2\vert = (j-1)(n-3),\ \ \ \vert E_3\vert = {j-1 \choose 2}.
 \end{equation}
Since $\vert E_1\vert = {n-2 \choose 2}$, using the proof of part (i) of Claim \ref{Ex} and  Lemma \ref{sdr}, for some $1 \leq i' \leq n-1$, $A_{i'} = \emptyset$ and
$$A_{i} = \{ e \in E(\mathcal{H}) : \  e=\{x,v_{i},v_l\}, \ \ l\neq i,i'\}, \ \ \ \ \ \  1\leq i \leq n-1\ \  {\rm and} \ \ i\neq i'.$$
 If $j \leq i' \leq n-1$, using (\ref{ddd}), we have  $\mathcal{H} \cong T_3(N,n-1)$. Hence we may assume that for some $1 \leq i' \leq j-1$, say $i'=1$, ${A}_1 = \emptyset$. By considering the sets $E_1$ and $E_2$ and using (\ref{ddd}), it can be shown that $\mathcal{H}[x,x_1,v_2,...,v_{n-1}]\cong \mathcal{K}_{n}^{3}$ and so it contains a  copy of   $K_{n}^{(3)}$. This contradiction completes the proof of the  theorem.

Now we may assume that for every vertex $x \in V(\mathcal{H})$, $d_{\mathcal{H}}(x) \geq d+1$. Set $G=B(\mathcal{H})$. So we may assume that $G=[X,Y]$, where  $$X=\{u_{ik}=v_iv_k :\ i< k\ \ {\rm and}\  \  v_i,v_k \in V(\mathcal{H})\}$$ and $Y=E(\mathcal{H})$.  Since, by (\ref{lbound}),  $|Y|\geq \sum _{i=0}^{3} \ell ^{3-i} {j \choose i}{n-1-i \choose 3-i},$ we have $\vert X\vert \leq \vert Y\vert$. Recall that every matching of $X$ in $G$ is equivalent to a ${K}_{N}^{(3)}$ in $\mathcal{H}$.  Let $X=X_1 \cup X_2 \cup X_3$, where
\begin{eqnarray*}
X_{1}&=&\{u \in X : d_{G}(u)=0\},\\
X_{2}&=&\{u \in X : 1 \leq d_{G}(u)\leq 4\},\\
X_{3}&=&\{u \in X : d_{G}(u)\geq 5\}.
\end{eqnarray*}
We have one of the following cases:\\

\noindent{\bf Case 1.} $X_1 \cup X_2 =\emptyset.$\\
In this case, the Hall's theorem \cite{hall} guarantees the existence of a matching of $X$ in $G$. That is a contradiction.\\\\

\noindent{\bf Case 2.}  $X_1 \cup X_2 \neq \emptyset.$\\
Let $X_1 \cup X_2 = \{v_{i_1}v_{i'_1},v_{i_2}v_{i'_2},...,v_{i_t}v_{i'_t}\}$. First we show that the following claim holds.

\begin{emp}\label{pairdisj2}
The elements of $X_1 \cup X_2$  are pairwise disjoint.
\end{emp}
\noindent\textbf{Proof of Claim \ref{pairdisj2}}. Suppose to contrary that for some  $2 \leq s \leq t$, $ \{wv_{i'_1},wv_{i'_2},...,wv_{i'_s}\} \subseteq X_1 \cup X_2$. So we have  $ d_{\mathcal{H}}(w) \leq f(s),$ where $f(k)= 4k + {n+j-k-2 \choose 2}$ is a function on $k$,  $2 \leq k \leq t \leq N-1$. It is straightforward to see that the absolute maximum of $f(k)$ occurs in point $k=2$. Hence $$ d_{\mathcal{H}}(w) \leq f(2) = 8 + {n+j-4 \choose 2}.$$
On the other hand,  $ 8 + {n+j-4 \choose 2} < d+1$ for $n\geq 13$. That  contradiction completes the proof of our claim.
$\hfill\square$\\

 Also we have the following claim.
 \begin{emp}\label{sizex2}
$\vert X_1 \vert \leq j-1$.
\end{emp}
\noindent\textbf{Proof of Claim \ref{sizex2}}.  Suppose not. Therefore we may assume that $\{v_{i_1}v_{i'_1},v_{i_2}v_{i'_2},...,v_{i_j}v_{i'_j}\}\subseteq X_1.$ Set $L=\{v_{i_2},v_{i_3},...,v_{i_j}\}$.  We have $E_{v_{i_1}}= F_1 \cup F_2 \cup F_3$, where
$$F_k =\{e \in E_{v_{i_1}} : \vert e \cap L\vert =k-1\},\ \ \ \ \ \ 1 \leq k \leq 3.$$
Since, using Claim \ref{pairdisj2}, the elements of $X_1 $  are pairwise disjoint, the elements of $L$ are distinct. So, an easy computation shows that     $\vert F_1\vert \leq {n-2 \choose 2}$, $\vert F_2\vert \leq (j-1)(n-3)$ and $\vert F_3\vert \leq {j-1 \choose 2}$. Therefore $$d_{\mathcal{H}}(v_{i_1})= \vert E_{v_{i_1}}\vert \leq {n-2 \choose 2}+(j-1)(n-3)+{j-1 \choose 2}=d.$$ This contradiction completes the proof of this claim.
$\hfill\square$\\

Using the definition of $X_2,$ for every $u_{ik}=v_iv_k\in X_2,$ we have $d_{G}(u_{ik}) \geq 1.$
 On the other hand, by  Claim \ref{pairdisj2}, the elements of $X_2$ are pairwise disjoint. Therefore  $G$ contains a matching $M_1$ of $X_2$ in $G$.
 Suppose that $G'=[X_3,Y']$ is the induced  subgraph of $G$  so that $Y'\subseteq Y$ is obtained by deleting the vertices of $M_1$. Note that for every $u \in X_3$ and $y \in Y'$, we have $d_{G'}(u) \geq 3$ and  $d_{G'}(y) \leq 3$. So the Hall's theorem \cite{hall} guarantees  the existence of a matching $M_2$ of $X_3$ in $G'$. Now, using Claim \ref{sizex2}, we may suppose that $X_1= \{v_{i_1}v_{i'_1},v_{i_2}v_{i'_2},...,v_{i_t}v_{i'_t}\}$, where $t\leq j-1$. Set $V'= V(\mathcal{H})\setminus \{v_{i'_1},v_{i'_2},...,v_{i'_t}\}$. Clearly
 $\vert V'| \geq n$ and   $M_1 \cup M_2$ induces a matching of $X_2\cup X_3$ in $G$. As every matching of $X_2\cup X_3$ in $G$ is equivalent to a
 $K_{\vert V' \vert}^{(3)}$ in $\mathcal{H}[V']$, we have a copy of  $K_{n}^{(3)}$ in $\mathcal{H}$. This is a contradiction to our assumption.


$\hfill\blacksquare$\\

\section{proof of Theorem \ref{main}}
Let $\mathcal{H}$ be an $\mathcal{F}_n^{(3)}$-free hypergraph with $N$ vertices and $ex(N,\mathcal{F}_n^{(3)})$ edges. Also let $N=\ell (n-1)+j,$ where $\ell\geq 1$ and $1\leq j \leq n-1.$
We use induction on $\ell$ to show that $ex(N,\mathcal{F}_n^{(3)})=t_3(N,n-1)$. Using Theorem \ref{l=1}, the basic step $\ell =1$ is true. Now suppose that $\ell >1$. Since at least one $K_{n}^{(3)}$ is made by adding one edge to $\mathcal{H}$, we deduce that $\mathcal{H}$ contains a $K_{n-1}^{(3)}$. Let $\mathcal{K}$ be such a $K_{n-1}^{(3)}$ in $\mathcal{H}$ with the core sequence $v_1,v_2,...,v_{n-1}$ so that $e(\mathcal{H}[v_1,v_2,...,v_{n-1}])\cap e(\mathcal{K})$ is maximum.
Let $\mathcal{H}_1=\mathcal{H}[V_1]$, $\mathcal{H}_2=\mathcal{H}[V_2]$ and $\mathcal{H}_3=\mathcal{H}[V_1 ,V_2]$, where $V_1 =V(\mathcal{K})=\{v_1,v_2,...,v_{n-1}\}$ and $V_2=V(\mathcal{H})\setminus V_1$. Also let $N'= \vert V_2\vert = (\ell -1)(n-1)+j$, where $\ell >1$ and $1\leq j \leq n-1$. Set
$$\mathcal{H}_3 ^{\vartriangle}= \{e \in E(\mathcal{H}_3) : \ \vert e \cap V_1 \vert =1\ \  {\rm and} \ \ \vert e \cap V_2 \vert =2\},$$
$$\mathcal{H}_3 ^{\triangledown}= \{e \in E(\mathcal{H}_3) :\  \vert e \cap V_1 \vert =2 \ \ {\rm and} \ \  \vert e \cap V_2 \vert =1\}.$$
Note that $E(\mathcal{H}_3)=\mathcal{H}_3 ^{\vartriangle} \cup  \mathcal{H}_3 ^{\triangledown}$.
So
\begin{equation}\label{eh}
e(\mathcal{H}) =e(\mathcal{H}_1) +e(\mathcal{H}_2) + \vert \mathcal{H}_3 ^{\vartriangle}\vert + \vert  \mathcal{H}_3 ^{\triangledown}\vert.
\end{equation}
By the induction hypothesis, we have
\begin{equation}\label{eh2}
e(\mathcal{H}_2) \leq t_3(N',n-1)=\sum _{i=0}^{3} (\ell-1) ^{3-i} {j \choose i}{n-1-i \choose 3-i}.
\end{equation}
Moreover,
\begin{equation}\label{eh3}
\vert \mathcal{H}_3^{\vartriangle} \vert \leq  t_2(N',n-1).
\end{equation}
To see that, let $G$ be a graph on $V_2$ so that the vertices $u$ and $v$ of $V_2$ are adjacent in $G$ if and only if there exists the edge $\{x,u,v\} \in \mathcal{H}_3 ^{\vartriangle}$, for some $x \in V_1$. If there is a $K_n$ in $G$, then we can find a $K_n^{(3)}$ in $\mathcal{H}$, a contradiction. Therefore, by Tur\'{a}n's theorem \cite{turan}, we have $\vert \mathcal{H}_3^{\vartriangle} \vert \leq  t_2(N',n-1)$.\\

\noindent Now we show that $e(\mathcal{H}_1)+\vert \mathcal{H}_3 ^{\triangledown}\vert \leq {n-1 \choose 3}+N'{n-2 \choose 2}$.
For this purpose, set $$\mathcal{B}_1 ^{\triangledown}= \{ e \in \mathcal{H}_3 ^{\triangledown} :\  e \in E(\mathcal{K}) \}$$ and $\mathcal{B}_2 ^{\triangledown}=  \mathcal{H}_3 ^{\triangledown}\setminus \mathcal{B}_1 ^{\triangledown}$. Clearly, we have
\begin{eqnarray}\label{B2}
\vert \mathcal{B}_2 ^{\triangledown}\vert \leq N'{n-2 \choose 2}.
\end{eqnarray}
To see that, choose an arbitrary  vertex $u \in V_2$. Set $x=u$ and $U=V_1=\{v_1,v_2,...,v_{n-1}\}$ and  $\mathcal{A}_u=\{A_{1}^u,A_{2}^u,...,A^u_{n-1}\}$, where
$$A^u_{i} = \{ e \in \mathcal{B}_2 ^{\triangledown} :\   \{u,v_{i}\} \subset e\}.$$
If $\mathcal{A}_u$ contains an SDR, then $u,v_1,v_2,...,v_{n-1}$ is the core sequence of a copy of  $K_n ^{(3)}$ in $\mathcal{H}$, a contradiction. So using Lemma \ref{sdr}, we have $\vert \bigcup _{i=1}^{n-1} A^u_{i} \vert \leq {n-2 \choose 2}$. Since $u$ is choosed as an arbitrary vertex of $V_2$, Thus $\vert \mathcal{B}_2 ^{\triangledown}\vert \leq N'{n-2 \choose 2} $. Now we demonstrate that 

\begin{eqnarray}\label{eh1h3}
e(\mathcal{H}_1) +\vert \mathcal{B}_1 ^{\triangledown}\vert \leq {n-1 \choose 3}.
\end{eqnarray}

To see this,  Suppose that $\vert \mathcal{B}_1 ^{\triangledown}\vert =t$. If  $t\leq e(\mathcal{H}_1 ^c)$, then we are done. So we may assume that $e(\mathcal{H}_1 ^c) \leq t-1$. On the other hand, clearly $\mathcal{K}_{n-1}^3$ contains a copy of  ${K}_{n-1}^{(3)}$. Therefore, by the maximality of $\mathcal{K}$, at most $t-1$ edges of $\mathcal{K}$ are not in $E(\mathcal{H}_1)$. This is a contradiction to the assumption that $\vert \mathcal{B}_1 ^{\triangledown}\vert =t$.

%
%

 Therefore by (\ref{B2}) and  (\ref{eh1h3}), we have
 \begin{eqnarray}\label{H1H3}
 e(\mathcal{H}_1)+\vert \mathcal{H}_3 ^{\triangledown}\vert \leq {n-1 \choose 3}+N'{n-2 \choose 2}.
   \end{eqnarray}
Now set
$$B= {n-1 \choose 3}+N'{n-2 \choose 2} + t_3(N',n-1) +t_2(N',n-1).$$
Hence by (\ref{eh}),(\ref{eh2}),(\ref{eh3}) and  (\ref{H1H3}), we have
\begin{eqnarray*}
\vert E(\mathcal{H})\vert \leq B &=& {n-1 \choose 3}+((\ell -1)(n-1)+j){n-2 \choose 2} + (\ell -1)^{3}{n-1 \choose 3}
 + j(\ell -1)^{2}{n-2 \choose 2}\\
 &+& (\ell -1) (n-3){j \choose 2} +{j \choose 3} + (\ell -1)^{2}{n-1 \choose 2}
+ j(\ell -1)(n-2)
+ {j \choose 2}.
\end{eqnarray*}

\noindent To demonstrate that $ex(N, \mathcal{F}_n^{(3)}) \leq t_3(N,n-1)$, it suffices to show that
$$B \leq t_3(N,n-1) =  \ell^{3}{n-1 \choose 3} + j\ell^{2}{n-2 \choose 2} + \ell (n-3){j \choose 2} +{j \choose 3}.$$
By simplifying the above inequality, it suffices to show that
\begin{eqnarray*}
3\ell {n-1 \choose 3} + (\ell -1)^2{n-1 \choose 2}+j(\ell -1)(n-2)\\ + {j \choose 2} +(\ell -1)(n-1){n-2 \choose 2} +2j{n-2 \choose 2}\\
 \leq 3 \ell^{2} {n-1 \choose 3}+2j\ell {n-2 \choose 2}+{j \choose 2}(n-3).
\end{eqnarray*}
But the above inequality is certainly true since $n \geq 13$, $\ell >1$ and $j \geq 1$ imply
\begin{itemize}
\item[$\bullet$]
${j \choose 2} \leq {j \choose 2}(n-3)$.
\item[$\bullet$]
 $ j(\ell -1)(n-2) +2j{n-2 \choose 2} \leq 2j\ell {n-2 \choose 2}$.
\item[$\bullet$]
$3\ell {n-1 \choose 3}+ (\ell -1)^2{n-1 \choose 2}+(\ell -1)(n-1){n-2 \choose 2} \leq 3 \ell^{2} {n-1 \choose 3}$.
\end{itemize}
So $$ex(N, \mathcal{F}_n^{(3)}) \leq t_3(N,n-1)$$ and the equality follows by inspection of (\ref{lbound}). Therefore,
\begin{itemize}
\item[(i)] $e(\mathcal{H}_{2}) =t_3(N',n-1)$.
\item[(ii)] $\vert \mathcal{H}_3^{\vartriangle} \vert =  t_2(N',n-1)$.
\item[(iii)] $\vert \mathcal{B}_2 ^{\triangledown}\vert = N'{n-2 \choose 2}$.
\item[(iv)] $e(\mathcal{H}_1) +\vert \mathcal{B}_1 ^{\triangledown}\vert = {n-1 \choose 3}$.
\end{itemize}

In the sequel, we demonstrate that $\mathcal{H}\cong T_3(N,n-1)$. Since  $e(\mathcal{H}_{2}) =t_3(N',n-1)$,  the induction hypothesis implies that $\mathcal{H} _{2}\cong T_3(N',n-1)$. Therefore $\mathcal{H}_2$ is a complete $3$-uniform $(n-1)$-partite hypergraph  on $N'$ vertices whose partition sets differ in size by at most 1.
 Assume that $U_1,U_2,...,U_{n-1}$  are the partition sets of $\mathcal{H}_2$. Recall that  $U=V_1=\{v_1,v_2,...,v_{n-1}\}$ and  for every $u \in V_2$,  $\mathcal{A}_u=\{A_{1}^u,A_{2}^u,...,A^u_{n-1}\}$, where
$A^u_{i} = \{ e \in \mathcal{B}_2 ^{\triangledown} :\   \{u,v_{i}\} \subset e\}$ and $\vert \bigcup _{i=1}^{n-1} A^u_{i} \vert \leq {n-2 \choose 2}$.  
  Since  $\vert \mathcal{B}_2 ^{\triangledown}\vert = N'{n-2 \choose 2}$,
  $$\vert \bigcup _{i=1}^{n-1} A^{u}_{i} \vert = {n-2 \choose 2},\ \ \ \ \ \ \forall \  u \in V_2.$$
   So  using Lemma \ref{sdr}, there exists $1 \leq q_u \leq n-1$, so that $A^{u}_{q_{u}}=\emptyset$ and for every  $1 \leq i \leq n-1$ and $i \neq q_u$, we have
$$A_{i}^{u} = \{ e \in \mathcal{B}_2 ^{\triangledown} :\  \{u,v_{i}\} \subset e\}=\{  \{u,v_i,v_k\}:\ k \neq i, q_u \}.$$
 
 On the other words,$$\bigcup _{i\neq q_u} A_i^{u}= \{ \{u,v_l,v_k\}:\ v_l ,v_k \in V_1,\ l,k \neq q_u\}.$$
 So we can partition the vertices of $V_2$ into $n-1$ parts $U'_1,U'_2,...,U'_{n-1}$, so that for every $x \in U'_m$, $A^{x}_m = \emptyset$ and
 $$\bigcup _{i\neq m} A_i^{x}= \{ \{x,v_l,v_k\}:\ v_l ,v_k \in V_1,\ l,k \neq m\}.$$
 Now we show that for $1\leq i\leq n-1$, $U'_i$ is an independent set in $\mathcal{H}$. Suppose not. By symmetry we may assume that for two vertices $x,y \in U'_1$, the edge $\{x,y,z\} \in E(\mathcal{H})$. It can be shown that $x,y,v_2,v_3,...,v_{n-1}$ represents the core sequence of a copy of  $K_n^{(3)}$ in $\mathcal{H}$ with the following edge assignments. Set $e_{xy}=\{x,y,z\}$, $e_{xv_i}\in A_i^x$ for $2\leq i\leq n-1$, $e_{yv_i}\in {A}_i^y$ for $2\leq i\leq n-1$ and $e_{v_iv_{i'}}\in E(\mathcal{K})$ for $2\leq i,i'\leq n-1$. Hence $U'_i$'s, $1 \leq i \leq n-1$, are independent sets in $\mathcal{H}$.\\

 Therefore $\{U_1,U_2,...,U_{n-1}\}=\{U'_1,U'_2,...,U'_{n-1}\}$. With no loss of generality, we may suppose that
 $$U_i =U'_i \ \ \ \ \ \ \ \ {\rm for}\ \  1\leq i\leq n-1.$$
 Now we demonstrate that for $1\leq i\leq n-1$, $U_i\cup \{v_i\}$ is an independent set in $\mathcal{H}$. Suppose to the contrary that for some $1\leq h\leq n-1$, $U_h \cup \{v_h\}$ is not independent set. So for some $u_h \in U_h$, $f=\{u_h,v_h,w\}\in E(\mathcal{H})$. Since $U_h$ is an independent set in $\mathcal{H}$, $w\notin U_h$.  Choose the vertices $x_1,x_2,...,x_{n-1}$ so that $x_h=u_h$ and
 $$x_i \in U_i\ \ \ \ \ \ \ \  {\rm for} \ \ i\neq h.$$
 Since $\mathcal{H}[\{x_1,x_2,...,x_{n-1}\}]\cong \mathcal{K}_{n-1}^3$, $x_1,x_2,...,x_{n-1}$ is the core sequence of a ${K}_{n-1}^{(3)}$, say $\mathcal{K}'$, in $\mathcal{H}$. Thus $x_1,x_2,...,x_{n-1},v_h$ represents the core sequence of a ${K}_n^{(3)}$ in $\mathcal{H}$ with the following edge assignments. Set $e_{x_hv_h}=f$, $e_{v_hx_i}\in A_h^{x_i}$ for $i\neq h$ and $e_{x_i x_{i'}}\in E(\mathcal{K}')$ for $i,i' \neq h$.

\noindent Therefore for $1\leq i \leq n-1$, $W_i=U_i \cup \{v_i\}$ is an independent set in $\mathcal{H}$. So   $\mathcal{H}$ is an $n-1$-partite hypergraph with parts $W_1,W_2,...,W_{n-1}$  whose partition sets differ in size by at most 1. Since $e(\mathcal{H})=t_3(N,n-1)$, we deduce that $\mathcal{H}\cong T_3(N,n-1)$.


$\hfill\blacksquare$\\



\footnotesize
\begin{thebibliography}{99}



\bibitem{Erdos-Gallai}  P. Erd{\H{o}}s, T. Gallai, On maximal paths and circuits of graphs, {\it Acta Math. Acad. Sci. Hungar.} {\bf 10} (1959) 337-356.

\bibitem{Erdos-ko-rado} P. Erd{\H{o}}s,  C. Ko and  R. Rado, Intersection theorems for systems of finite sets, \textit{Quart. J. Math. Oxford Ser. (2)} {\bf 12}  (1961), 313--320.

\bibitem{loose pathI} Z. F{\"u}redi,  T. Jiang and R. Seiver, On Exact solution of the hypergraph {T}ur\'an problem for {$k$}-uniform linear paths,\textit{ Combinatorica} {\bf 34} (2014), 299--322.

\bibitem{loose cycleI} Z. F{\"u}redi and   T. Jiang, Hypergraph {T}ur\'an numbers of linear cycles, \textit{J. Combin. Theory Ser. A} {\bf 123} (2014), 252--270.

\bibitem{Berge cycles} Z. F{\"u}redi and L. {\"O}zkahya, On 3-uniform hypergraphs without a cycle of a given length,  \textit{ arXiv:1412.8083v2}.


\bibitem{Berge paths} E. Gy{\H{o}}ri, G. K. Katona  and N. Lemons, Hypergraph extensions of the {E}rd{\H o}s-{G}allai theorem, \textit{European J. Combin.} 58 (2016) 238--246.


\bibitem{Berge cycles2} E. Gy{\H{o}}ri, N. Lemons, Hypergraphs with no cycle of a given length, \textit{ Combin. Probab. Comput.} {\bf 21} (2012),
193--201.

\bibitem{Berge cyclesI} E. Gy{\H{o}}ri and N. Lemons, 3-uniform hypergraphs avoiding a given odd cycle,\textit{ Combinatorica} {\bf 32} (2012), 187--203.



\bibitem{hall} P. Hall, On representatives of subsets, \textit{J. London Math. Soc.} {\bf 10} (1935), 26--30.
\bibitem{survey} P. Keevash, Hypergraph Tur\'{a}n problems, in: \textit{Surveys in Combinatorics 2011}, Cambridge
University Press, (2011), 83--140.

\bibitem{loose paths and cycles} A. Kostochka, D. Mubayi and J. Verstra{\"e}te, Tur\'an problems and shadows {I}: {P}aths and cycles, \textit{J. Combin. Theory Ser. A}
{\bf 129} (2015), 167--190.

\bibitem{mubay} D. Mubayi, A hypergraph extension of Tur\'{a}n's theorem,  \textit{ J. Combin. Theory
Ser. B} {\bf 96} (2006), 122--134.

\bibitem{pikh} O. Pikhurko, Exact computation of the hypergraph Turán function for expanded complete 2-graphs, \textit{ J. Combin.
Theory Ser. B}, {\bf 103} (2013), 220--225.

\bibitem{turan} P. Tur{\'a}n, Eine {E}xtremalaufgabe aus der {G}raphentheorie, \textit{Mat. Fiz. Lapok} {\bf 48} (1941) 436--452.



\end{thebibliography}
\end{document}